\documentclass[11pt,letterpaper]{amsart}

\usepackage[centertags]{amsmath}
\usepackage{amsthm}
\usepackage{amssymb}
\usepackage{amscd}
\usepackage{eucal}
\usepackage[matrix,arrow]{xy}

\newcommand{\R}{\mathbb{R}}
\newcommand{\Z}{\mathbb{Z}}

\newcommand{\iso}{\cong}           
\newcommand{\htp}{\simeq}          
\newcommand{\smooth}{C^\infty}
\newcommand{\suchthat}{\; | \;}

\newcommand{\id}{\mathrm{id}}
\newcommand{\im}{\mathrm{im}}
\renewcommand{\o}{\omega}
\newcommand{\Symp}{\mathrm{Symp}}
\newcommand{\Diff}{\mathrm{Diff}}

\theoremstyle{plain}
\newtheorem{theorem}{Theorem}
\newtheorem{lemma}[theorem]{Lemma}
\newtheorem{definition}[theorem]{Definition}

\newcommand{\Sympm}{\Symp^m}
\newcommand{\tphi}{\tilde{\phi}}
\newcommand{\gen}[1]{\langle #1 \rangle}
\newcommand{\tx}{\tilde{x}}
\newcommand{\ZZ}{\mathcal Z}
\newcommand{\JJ}{\mathcal J}
\newcommand{\BB}{\mathcal B}
\newcommand{\MM}{\mathcal M}
\newcommand{\duds}{\frac{\partial u}{\partial s}}
\newcommand{\dudt}{\frac{\partial u}{\partial t}}
\newcommand{\tM}{\tilde{M}}

\title[Floer homology and mapping class group]{Symplectic
Floer homology\\ and the mapping class group}
\author{Paul Seidel}
\date{Revised version, March 3, 2001.}

\begin{document}
\begin{abstract}
We consider symplectic Floer homology in the lowest nontrivial dimension, that is
to say, for area-preserving diffeomorphisms of surfaces. Particular attention is paid
to the quantum cap product.
\end{abstract}
\maketitle

\section{Introduction}

Let $\Gamma = \pi_0(\Diff^+(M))$ be the mapping class group of a closed connected
oriented surface $M$ of genus $\geq 2$. Pick an everywhere positive two-form $\o$
on $M$. A theorem of Moser says that each $g \in \Gamma$ admits representatives
which preserve $\o$. These are automorphisms of the two-dimensional symplectic
manifold $(M,\o)$, and the techniques of symplectic topology can be applied to
them. We will use symplectic representatives $\phi$ which satisfy
a certain additional monotonicity property. The benefit is that the Floer homology
$HF_*(\phi)$, which is a finite-dimensional $\Z/2$-graded vector space over the field
$\Z/2$, is independent of the choice of $\phi$, hence an invariant of $g$; we denote it by
$HF_*(g)$. There is an additional multiplicative structure, the quantum cap product
\[
\ast: H^*(M;\Z/2) \otimes HF_*(\phi) \rightarrow HF_*(\phi);
\]
this is again independent of the choice of representative $\phi$, thus
equipping $HF_*(g)$ with the structure of a $\Z/2$-graded $H^*(M;\Z/2)$-module.
The most familiar case is when $g = e$: $HF_*(e) \iso H_*(M;\Z/2)$, and $\ast$
reduces to the ordinary cap product (this was announced in
\cite{ruan-tian-bott} and \cite{piunikhin-salamon-schwarz94}; a detailed proof,
in greater generality than envisaged by the earlier approaches, has appeared in
\cite{liu-tian-final}). In particular, the quantum cap action of $H^2(M;\Z/2)
\iso \Z/2$ on $HF_*(e)$ is nonzero. Our first result is that this characterizes
the trivial mapping class.

\begin{theorem} \label{th:one}
For all $g \neq e$ in $\Gamma$, the
quantum cap action of $H^2(M;\Z/2)$ on $HF_*(g)$ is zero.
\end{theorem} \vspace{-0.5em}

Floer homology does not distinguish between arbitrary mapping classes:
there are cases when $HF_*(g_1) \iso HF_*(g_2)$ as $H^*(M;\Z/2)$-modules and
where $g_1,g_2$ are not even conjugate. An easy source of such examples are fixed
point free diffeomorphisms of finite order, since their Floer homology
is always zero.

\begin{theorem} \label{th:two}
Suppose that $a \in H^1(M;\Z/2)$ is a class whose quantum cap action on
$HF_*(g)$ is nonzero. Then there is an $l: S^1 \rightarrow M$ with $g(l) \htp
l$ and $\langle a, [l] \rangle = 1$.
\end{theorem} \vspace{-0.5em}

Here $\htp$ denotes the relation of free homotopy between loops in $M$. The same
examples as before can be used to see that the converse to Theorem \ref{th:two} is false.
By this we mean that one can have $g(l) \htp l$ with $[l] \neq 0$ in $H_1(M;\Z/2)$, and the
action of $H^1(M;\Z/2)$ on $HF_*(g)$ can still be zero.

The strategy of proof is to consider a certain infinite (not usually Galois or
connected) covering $p_\phi: \tM_\phi \rightarrow M$ canonically associated to
$\phi \in \Diff^+(M)$. One can introduce a refined quantum cap product, in which
the cohomology of this covering replaces that of $M$. More precisely:

\begin{lemma} \label{th:viterbo}
Let $\phi$ be a monotone symplectic automorphism of $(M,\o)$. Then there is
a map $\tilde{\ast}$ fitting into a commutative diagram
\[
\xymatrix{
 {H^*(\tM_\phi;\Z/2) \otimes HF_*(\phi)} \ar[rr]^-{\tilde{\ast}} &&
 {HF_*(\phi)} \\
 {H^*(M;\Z/2) \otimes HF_*(\phi)} \ar[u]^{p_\phi^* \otimes \id} \ar[rru]^{\ast}
}
\]
In particular, if the quantum cap action of $a \in H^*(M;\Z/2)$ is nonzero, $p_\phi^*(a)$
must be nonzero.
\end{lemma}

This, together with elementary topological arguments, yields the
theorems above. To round off the discussion, we should mention that Floer homology
can also be defined with integer coefficients. Our arguments carry over
without any changes, and in the case of Theorem \ref{th:two} this leads to a
slightly better control over the homology class of $l$.

{\em Acknowledgments.}
The importance of the quantum module structure for symplectic isotopy questions
was first pointed out by Callahan. Lemma \ref{th:viterbo} is inspired by a
construction of Viterbo
\cite{viterbo95}. Ivan Smith and Dietmar Salamon made useful comments on a
preliminary version of this paper.

\section{Three topological preliminaries\label{sec:topo}}

(a) Let $\phi \in \Diff^+(M)$. The definition of the covering $p_\phi: \tM_\phi
\rightarrow M$ is that points of $\tM_\phi$ are pairs $(x,[c])$, where $x
\in M$, $c: [0;1] \rightarrow M$ is a path from $c(0) = \phi(x)$ to $c(1) = x$,
and $[c]$ is the homotopy class of $c$ rel endpoints. Equivalently, after
choosing a hyperbolic metric on $M$, one can describe $\tM_\phi$ as the space
of all geodesic paths $c: [0;1] \rightarrow M$ such that $c(0) = \phi(c(1))$.
There is a canonical lift of $\phi$ to a diffeomorphism $\tphi: \tM_\phi
\rightarrow \tM_\phi$, $\tphi(x,[c]) = (\phi(x), [\phi(c)])$; and it is not
difficult to see that $\tphi$ is homotopic to the
identity map. In particular, if a loop $l: S^1 \rightarrow M$ can be lifted to
$\tM_\phi$ then $\phi(l) \htp l$; and the converse is also true. This has an
immediate consequence:

\begin{lemma} \label{th:pullback1}
Let $a \in H^1(M;\Z/2)$ be a class such that $p_\phi^*(a) \neq 0$. Then there
is a loop $l$ such that $\phi(l) \htp l$ and $\langle a,[l] \rangle = 1$. \qed
\end{lemma} \vspace{-0.5em}

Now suppose that $\tM_\phi$ has a connected component which is compact, so that
the restriction of $p_\phi$ to it is a covering of finite order.  Then every
loop $l$ has an iterate $l^k$, $k \geq 1$, which lifts to $\tM_\phi$, so that
$\phi(l)^k \htp l^k$. As one can see by looking at the unique geodesic
representatives, this implies $\phi(l) \htp l$. From that one can deduce in
various ways that $[\phi] \in \Gamma$ is trivial, for instance: since lengths
of closed geodesics form a set of coordinates on the Teichm{\"u}ller space
${\mathcal T}$ of $M$ \cite[Expos{\'e} 7, Th{\'e}or{\`e}me
4]{fathi-laudenbach-poenaru}, it follows that the class $g = [\phi] \in \Gamma$
acts trivially on ${\mathcal T}$. In other words, every hyperbolic metric on
$M$ admits an isometry lying in $g$. This means that $g = e$, or else that $M$
is of genus two and $\phi$ is isotopic to a hyperelliptic involution. However,
hyperelliptic involutions act as $-1$ on $H_1(M;\Z)$, hence they can not
satisfy $\phi(l) \htp l$ for homologically nontrivial $l$. We have therefore
proved:

\begin{lemma} \label{th:pullback2}
If $H^2(\tM_\phi;\Z/2)$ is nonzero, $[\phi] \in \Gamma$ is the identity class.
\qed
\end{lemma} \vspace{-0.5em}

(b) For $\phi$ in the group $\Symp(M,\o)$ of symplectic automorphisms, consider
its mapping torus $T_\phi = (\R \times M) / (t+1,x) \sim (t,\phi(x))$, which is
a three-manifold fibered over $S^1$. The pullback of $\o$ to $\R \times M$
descends to a closed two-form $\o_\phi \in \Omega^2(T_\phi)$. The tangent
bundle along the fibres of $T_\phi \rightarrow S^1$ is an oriented real
two-plane bundle, whose Euler class we denote by $c_\phi \in H^2(T_\phi;\R)$.
Now $H^*(T_\phi;\R)$ fits into an exact sequence
\begin{equation} \label{eq:torus}
\cdots \rightarrow H^1(M;\R) \xrightarrow{\id - \phi^*} H^1(M;\R) \stackrel{d}
\rightarrow H^2(T_\phi;\R) \twoheadrightarrow H^2(M;\R),
\end{equation}
where the last map is restriction to the fibre. Since the difference $[\o_\phi]
- (\int_M \o/\chi(M)) \, c_\phi$ vanishes when restricted to a fibre, it is of
the form $d(m(\phi))$ for some unique $m(\phi) \in H^1(M;\R)/\im(\id -
\phi^*)$. A computation shows the following:

\begin{lemma} \label{th:calabi}
Let $(\psi_t)_{0 \leq t \leq 1}$ be an isotopy in $\Symp(M,\o)$ with $\psi_0 =
\id$, and $c \in H^1(M;\R)$ its Calabi invariant. Then for any $\phi \in
\Symp(M,\o)$ one has $m(\phi \circ \psi_1) = m(\phi) + c$. \qed
\end{lemma} \vspace{-0.5em}

Call $\phi$ monotone if $m(\phi) = 0$. These maps form a closed subspace (not
subgroup) $\Sympm(M,\o) \subset \Symp(M,\o)$, and using Lemma \ref{th:calabi}
one can show that this is a deformation retract of $\Symp(M,\o)$. By combining
this with Moser's theorem \cite{moser65} that $\Symp(M,\o) \hookrightarrow
\Diff^+(M)$ is a homotopy equivalence, and the Earle-Eells theorem
\cite{earle-eells67} that the connected components of $\Diff^+(M)$ are
contractible, one gets:

\begin{lemma} \label{th:moser}
The inclusion $\Sympm(M,\o) \hookrightarrow \Diff^+(M)$ induces an isomorphism
$\pi_0(\Sympm(M,\o)) \iso \Gamma$. Moreover, each connected component of
$\Sympm(M,\o)$ is contractible. \qed
\end{lemma} \vspace{-0.5em}

(c) Choose a Morse function $f$ on $M$, and a Riemannian metric such that the
gradient flow of $f$ is Morse-Smale. The Morse cochain complex $CM^*(f)$ is the
graded vector space over $\Z/2$ freely generated by the critical points of $f$,
with the degree of a generator $\gen{x}$ given by its Morse index $\mu(x)$. For
$\mu(x_+) = \mu(x_-)+1$, let $l(x_-,x_+) \in \Z/2$ be the number modulo 2 of
flow lines of $\nabla f$ going from $x_-$ to $x_+$. One defines a differential
$\delta_f: CM^*(f) \rightarrow CM^{*+1}(f)$ by $\delta_f\gen{x_-} = \sum_{x_+}
l(x_-,x_+) \gen{x_+}$. The cohomology $H(CM^*(f),\delta_f)$ is canonically
isomorphic to $H^*(M;\Z/2)$; see \cite{schwarz97} for a proof.

Now let $p: \tM \rightarrow M$ be a covering, and $\tilde{f} = f \circ p$.
Consider the graded vector space $CM^*(\tilde{f})$ whose elements are (possibly
infinite) formal sums of critical points of $\tilde{f}$, with
$\Z/2$ coefficients. Counting flow lines of $\nabla \tilde{f}$ with respect to
the pullback metric yields numbers $\tilde{l}(\tx_-,\tx_+) \in \Z/2$ for any
two critical points $\tx_-,\tx_+$ of $\tilde{f}$ with $\mu(\tx_+) =
\mu(\tx_-)+1$, and from those numbers one constructs a differential
$\delta_{\tilde{f}}$ as before. The argument from \cite{schwarz97} can be
adapted to show that there is a canonical isomorphism
\[
H(CM^*(\tilde{f}),\delta_{\tilde{f}}) \iso H^*(\tM;\Z/2).
\]
Since any flow line of $\nabla\tilde{f}$ is the lift of a flow line of $\nabla
f$, one finds that for any critical points $x_-,x_+$ of $f$ and any $\tx_+ \in
p^{-1}(x_+)$,
\[
l(x_-,x_+) = \sum_{p(\tx_-) = x_-} \tilde{l}(\tx_-,\tx_+).
\]
This implies that the homomorphism $p^*: CM^*(f) \rightarrow CM^*(\tilde{f})$
which takes any critical point to the sum of its preimages is a chain map. Via
the canonical isomorphisms, its induced map can be identified with the usual
pullback $p^*: H^*(M;\Z/2) \rightarrow H^*(\tM;\Z/2)$ on cohomology.

\section{Floer homology\label{sec:floer}}

Fix $\phi \in \Sympm(M,\o)$. Let $\Lambda_\phi = \{ y \in \smooth(\R,M)
\suchthat y(t) = \phi(y(t+1)) \}$ be the twisted free loop space, which is also
the space of sections of $T_\phi \rightarrow S^1$. The action form is the
closed one-form $\alpha_\phi$ on $\Lambda_\phi$ defined by
\[
\alpha_\phi(y) Y = \int_0^1 \o(dy/dt,Y(t))\,dt.
\]
The main difference between our exposition and the ones in the existing
literature, such as \cite{dostoglou-salamon94}, is that we will admit a
larger-than-usual class of perturbations of this one-form.

\begin{definition}
$\BB$ is the space of smooth families $b = (b_t)_{t \in \R}$ of closed one-forms
on $M$ satisfying $b_{t+1} = \phi^*b_t$ and
\begin{equation} \label{eq:periods}
\int_0^1 [b_t] \, dt \in \im(\id-\phi^*) \subset H^1(M;\R).
\end{equation}
\end{definition}

A condition equivalent to \eqref{eq:periods} is that there should be closed
one-forms $a_t$ and functions $H_t$ on $M$ such that $a_{t+1} = \phi^*a_t$,
$H_{t+1} = \phi^*H_t$, and $b_t = -\partial a_t/\partial t + dH_t$. For a
more intrinsic formulation consider the closed two-form on $T_\phi$
whose pullback to $\R \times M$ is $dt \wedge b_t$. Then \eqref{eq:periods}
is equivalent to saying that this two-form is exact; in fact it is $-d(a_t
+ H_t dt)$. Define the $b$-perturbed action form to be
\[
\alpha_{\phi,b}(y) Y = \alpha_\phi(y)Y + \int_0^1 b_t(Y(t)) \, dt.
\]
The perturbation term which we have added is exact, since it is the derivative of
the function $y \mapsto \int_0^1
a_t(dy/dt) + H_t(y(t))\,dt$ on $\Lambda_\phi$. Hence $[\alpha_{\phi,b}] = [\alpha_\phi] \in
H^1(\Lambda_\phi;\R)$. The set $\ZZ(b)$ of zeros of $\alpha_{\phi,b}$ consists
of the solutions $y \in \Lambda_\phi$ of
\begin{equation} \label{eq:critical}
dy/dt = X_{b,t}(y(t)),
\end{equation}
where $X_{b,t}$ is the symplectic vector field dual to $b_t$. There is an open
dense subset $\BB^{reg} \subset \BB$ such that if $b$ lies in it, all solutions
of \eqref{eq:critical} are nondegenerate. Nondegeneracy means that the map
$L_y: TM_{y(0)} \rightarrow TM_{y(1)}$ obtained by linearizing
\eqref{eq:critical} near $y$ satisfies $\det(\id - D\phi_{y(1)} \circ L_y) \neq
0$. One then defines $CF^*(\phi,b)$ to be the vector space over $\Z/2$ freely
generated by the finite set $\ZZ(b)$, with a $\Z/2$-grading given by the sign
of the above-mentioned determinant.

Let $\JJ$ be the space of smooth families $J = (J_t)_{t \in \R}$ of
$\o$-compatible almost complex structures on $M$ such that $J_{t+1} =
\phi^*J_t$. For $b \in \BB^{reg}$, $J \in \JJ$, and $y_-,y_+ \in \ZZ(b)$, let
$\MM(y_-,y_+,b,J)$ be the set of maps $u: \R^2 \rightarrow M$ which satisfy
\begin{equation} \label{eq:floer}
\left\{
\begin{split}
 & u(s,t) = \phi(u(s,t+1)), \\
 & \textstyle{\lim_{s \rightarrow \pm\infty}} u(s,\cdot) = y_{\pm}, \\
 & \partial u/\partial s + J_t(u)(\partial u/\partial t - X_{b,t}(u)) = 0.
\end{split}
\right.
\end{equation}
To each such map is associated a Fredholm operator $D_u$ which linearizes
\eqref{eq:floer} in suitable Sobolev spaces. We write $\MM_k(y_-,y_+,b,J)$ for
the subset of those $u$ with $\mathrm{index}(D_u) = k \in \Z$.

\begin{lemma} \label{th:energy}
The energy
\[
E(u) = \int_{\R \times [0;1]} \left| \duds \right|_{J_t}^2
\]
is constant on each $\MM_k(y_-,y_+,b,J)$.
\end{lemma} \vspace{-0.5em}

\proof Take $u,v \in \MM_k(y_-,y_+,b,J)$ and consider them as maps $\R
\rightarrow \Lambda_\phi$, $s \mapsto u(s,\cdot)$ respectively $v(s,\cdot)$,
with the same asymptotic behaviour. Because of the
gradient flow interpretation of \eqref{eq:floer}, the energy of each is equal
to minus the integral of $\alpha_{\phi,b}$ along it. Gluing $u,v$ together at
the ends yields a map $w = u \# \bar{v}: S^1 \rightarrow \Lambda_\phi$, unique
up to homotopy, such that $E(u)-E(v) = -\int_{S^1} w^*\alpha_{\phi,b} =
-\int_{S^1} w^*\alpha_\phi$. One may view $w$ as a map $S^1 \times S^1
\rightarrow T_\phi$, and rewrite the last equality as
\begin{equation} \label{eq:energy}
E(u)-E(v) = \int_{S^1 \times S^1} w^*\o_\phi.
\end{equation}
The fact that $\mathrm{index}(D_u) = \mathrm{index}(D_v) = k$, together with a
standard gluing theorem for the index, shows that $w^*c_\phi \in H^2(S^1 \times
S^1)$ is zero. Since $\phi$ is monotone this implies that the right hand side
of \eqref{eq:energy} is zero, so that $E(u) = E(v)$. \qed

For each $b \in \BB^{reg}$ there is a dense subset $\JJ^{reg}(b) \subset \JJ$
of almost complex structures $J$ such that all solutions of \eqref{eq:floer}
with arbitrary $y_-,y_+$ are regular, meaning that $D_u$ is onto. Then each
$\MM_k(y_-,y_+,b,J)$ has a natural structure of a smooth $k$-dimensional
manifold. Moreover, the quotients $\MM_1(y_-,y_+,b,J)/\R$ by translation in
$s$-direction are actually finite sets. Denoting by $n(y_-,y_+)$ the number of
points mod 2 in each of them, one defines a differential $\partial_{b,J}:
CF_*(\phi,b) \rightarrow CF_{* \pm 1}(\phi,b)$ by $\partial_{b,J}\gen{y_-} =
\sum_{y_+} n(y_-,y_+) \gen{y_+}$. This satisfies $\partial_{b,J} \circ
\partial_{b,J} = 0$, and its homology is by definition the Floer homology
$HF_*(\phi)$. The necessary analysis is the same as in the slightly different
situations considered in \cite{dostoglou-salamon94}, \cite{floer-hofer-salamon95}.
In fact, the pertubation does not really affect the analytic arguments,
since a change of variables transforms \eqref{eq:floer}
into a similar equation without the $X_{b,t}$ term. Note also that bubbling is not an
issue, as we have $\pi_2(M) = 0$.

A continuation argument shows that Floer homology is independent of the choice
of $b$ and $J$ up to canonical isomorphism. This involves an equation, with
periodicity and convergence conditions as before,
\begin{equation} \label{eq:continuation}
\duds + J_{s,t}(u)\big(\dudt - X_{b,s,t}(u)\big) = 0.
\end{equation}
Here the $b_{s,t}$ are closed one-forms on $M$, which are $s$-independent if
$|s|$ is sufficiently large; and one supposes again that $b_{s,t} =
-\partial a_{s,t}/ \partial t + dH_{s,t}$, where $a_{s,t}$ and $H_{s,t}$ are
also $s$-independent for $|s| \gg 0$. The construction follows closely the
familiar pattern \cite{salamon-zehnder92} and only one step deserves a separate
mention, namely, the proof of the energy estimate. Let $u$ be a solution of
\eqref{eq:continuation} with limits $y_{\pm}$. A straightforward computation yields
\begin{align*}
 & E(u) = \int_{\R \times [0;1]} \left| \duds \right|^2
 = \int_{\R \times [0;1]} \o\big(\duds,\dudt\big)
\\ &
 - \int_{\{+\infty\} \times [0;1]} \left( a_{s,t}\big(\dudt\big)
 + H_{s,t}(u)\right)
 + \int_{\{-\infty\} \times [0;1]} \left( a_{s,t}\big(\dudt\big) +
 H_{s,t}(u)\right)
\\ &
 + \int_{\R \times [0;1]} \frac{\partial H_{s,t}}{\partial s}(u)
 + \int_{\R \times [0;1]} \frac{\partial a_{s,t}}{\partial s}(X_{b,s,t}(u))
 + \int_{\R \times [0;1]} \frac{\partial a_{s,t}}{\partial s}
 \big(J_{s,t}(u)\duds\big).
\end{align*}
One can show in the same way as in Lemma \ref{th:energy} that the first
term $\int u^*\o$ depends only on $y_\pm$ and on the expected dimension of the
moduli space of solutions near $u$. The second and third term obviously depend
only on $y_\pm$; the fourth and fifth one are bounded by constants independent
of $u$; and the last one by a constant times $E(u)^{1/2}$. It follows that on
any relevant moduli space of solutions there is an inequality $E(u) \leq
{const}_1 + {const}_2\, E(u)^{1/2}$, giving an upper bound for $E(u)$.

Let $(\phi \circ \psi_t)_{0 \leq t \leq 1}$ be an isotopy in $\Sympm(M,\o)$
such that $\psi_t = \id_M$ for small $t$, and which is constant for $t$ near
$1$. Let $c_t$ be the closed one-forms on $M$ which generate the isotopy, $c_t =
\o(\cdot,D\psi_t^{-1}(\partial \psi_t/\partial t))$. Because of Lemma
\ref{th:calabi}, the fact that the isotopy stays within $\Sympm(M,\o)$ means
that
\begin{equation} \label{eq:restricted}
[c_t] \in \im(\id-\phi^*) \subset H^1(M;\R).
\end{equation}
Set $\phi' = \phi \circ \psi_1$ and define a diffeomorphism $\Psi:
\Lambda_{\phi'} \rightarrow \Lambda_\phi$ by $(\Psi y)(t) = \psi_t(y(t))$ for
$0 \leq t \leq 1$. A straightforward computation shows that
$\Psi^*\alpha_{\phi,b} = \alpha_{\phi',b'}$ where $b_t' = \psi_t^*b_t - c_t$
for $0 \leq t \leq 1$. \eqref{eq:restricted} implies that $b' = (b_t')
\in \BB$. Choosing $b \in \BB^{reg}$ and $J \in \JJ^{reg}(b)$ and setting $J'_t
= \psi_t^*J_t$, one has that $b' \in \BB^{reg}$ and $J' \in \JJ^{reg}(b')$.
$\Psi$ induces an isomorphism of chain complexes $(CF_*(\phi',b'),
\partial_{b',J'}) \rightarrow (CF_*(\phi,b),\partial_{b,J})$, which implies
that
\begin{equation} \label{eq:isotopy}
HF_*(\phi') \iso HF_*(\phi).
\end{equation}
This, together with the first part of Lemma \ref{th:moser}, allows one to
define $HF_*(g)$ for $g \in \Gamma$ as the Floer homology of any monotone
symplectic representative. One can actually prove that the isomorphisms
\eqref{eq:isotopy} are canonical, which means independent of the choice of the
path $(\phi \circ \psi_t)$ between $\phi$ and $\phi'$. We will not really use
this, but it seems appropriate at least to outline the argument. Suppose first
that we have an isotopy with $\psi_1 = \id_M$, so that $\phi' = \phi$. Then our
construction defines an automorphism of $HF^*(\phi)$. In a slightly different
situation, such automorphisms were studied in \cite{seidel96}, and a suitably
adapted version of \cite[Proposition 5.1]{seidel96} shows that they depend only
on the class of $(\psi_t)_{0 \leq t \leq 1}$ in $\pi_1(\Sympm(M,\o))$. This
fundamental group is trivial by Lemma \ref{th:moser}, hence the automorphism of
$HF^*(\phi)$ is the identity. The general statement can be derived from
this special case, by considering the composition of a map \eqref{eq:isotopy}
and the inverse of another one, coming from a different isotopy.

\section{Quantum cap product}

Take $\phi \in \Sympm(M,\o)$, $b \in \BB^{reg}$ and $J \in \JJ^{reg}(b)$ as
before. We begin by recalling
the Gromov-Floer compactification $\overline{\MM}_k(y_-,y_+,b,J)$
of the space $\MM_k(y_-,y_+,b,J)$. As a set, this is the disjoint union of
\begin{equation} \label{eq:broken}
\begin{split}
\prod_{i=1}^p \MM_{l_i}(y_-^{i-1},y_-^i,b,J)/\R \;\; \times & \;\;
\MM_{k-\sum_i l_i-\sum_j m_j}(y_-^p,y_+^0,b,J) \\
& \qquad \times
\prod_{j=1}^q \MM_{m_j}(y_+^{j-1},y_+^j,b,J)/\R
\end{split}
\end{equation}
for all $p,q \geq 0$, $l_i,m_j > 0$, and $y_-^i,y_+^j \in \ZZ(b)$ with $y_-^0 =
y_-$ and $y_+^q = y_+$. We will not write down the definition of the topology,
but one of its properties is that if a sequence $(u_k)$ in $\MM_k(y_-,y_+,b,J)$
has a limit $([v_1^-],\dots,[v_p^-],v,[v_1^+],\dots,[v_q^+])$ then the $u_k$
converge uniformly on compact subsets towards the principal component $v$
(which may be $s$-independent). As a consequence, the evaluation map $\eta:
\MM_k(y_-,y_+,b,J) \rightarrow M$, $\eta(u) = u(0,1)$, has a natural continuous
extension $\bar{\eta}$ to the compactification.

Choose a Morse function $f$ on $M$ and a Riemannian metric for which $\nabla f$ is
Morse-Smale. We may assume that the unstable (for the flow of $\nabla f$) manifolds
$W^u(f,x_-) \subset M$ of all critical points $x_-$ are transverse to $\eta$ for all
$y_-,y_+$ and $k$. If $k$ is equal to the Morse index $\mu(x_-)$, the preimage
$\eta^{-1}W^u(f,x_-) \subset \MM_k(y_-,y_+,b,J)$ is zero-dimensional. Moreover, a
dimension count using \eqref{eq:broken} shows that the two sets
\[
\left.
\begin{split}
\bar{\eta}^{-1}(\overline{W^u(f,x_-)} \setminus W^u(f,x_-)) \\
\bar{\eta}^{-1}W^u(f,x_-) \setminus \eta^{-1}W^u(f,x_-)
\end{split}
\right\}
\subset \overline{\MM}_k(y_-,y_+,b,J)
\]
are empty. This implies that $\eta^{-1}W^u(f,x_-) \subset \MM_k(y_-,y_+,b,J)$
is closed inside the Gromov-Floer compactification, hence a finite set. Let
$q(x_-,y_-,y_+)$ be the number of points mod 2 in it; one defines a
homomorphism
\begin{equation} \label{eq:chain}
\begin{split}
& CM^*(f) \otimes CF_*(\phi,b) \longrightarrow CF_*(\phi,b), \\ & \gen{x_-}
\otimes \gen{y_-} \longmapsto \sum_{y_+} q(x_-,y_-,y_+) \gen{y_+}.
\end{split}
\end{equation}
The next step is to consider the spaces $\eta^{-1}W^u(f,x_-) \subset \MM_k(y_-,y_+,b,J)$
in the case where they are one-dimensional, which is when $k = \mu(x_-) + 1$. A gluing
argument shows that their closures
$
\bar{\eta}^{-1}\overline{W^u(f,x_-)} \subset \overline{\MM}_k(y_-,y_+,b,J)
$
are one-dimensional manifolds with boundary, and by counting their boundary points one
sees that
\begin{multline*}
\sum_{x} l(x_-,x) q(x,y_-,y_+) + \sum_{y} n(y_-,y) q(x_-,y,y_+) \\ + \sum_{y}
n(y,y_+) q(x_-,y_-,y) = 0, \qquad \qquad
\end{multline*}
where $l$ and $n$ are defined in sections \ref{sec:topo}c and \ref{sec:floer}, respectively. This
means that \eqref{eq:chain} is a chain map. The induced map on homology is, by
definition, the quantum cap product $\ast: H^*(M;\Z/2) \otimes HF_*(\phi) \rightarrow HF_*(\phi)$.
This can be proved to be independent of the Morse function $f$ and the Riemannian
metric, and of the choices made in the definition of $HF_*(\phi)$. One can
also prove, but we will not really use it here, that it makes $HF_*(\phi)$ into a
module over the cohomology ring of $M$. Moreover the quantum cap product
commutes with the isomorphisms \eqref{eq:isotopy}, so that it gives rise to
a well-defined module structure on $HF_*(g)$ for $g \in \Gamma$. Details can
be found in \cite{le-ono95}, \cite{piunikhin-salamon-schwarz94},
\cite{schwarz97}, \cite{liu-tian-final}; the situation considered in those references
tends to differ slightly from the one here, but the arguments carry over with minimal
changes.

The construction of the refined product $\tilde{\ast}$
parallels that of $\ast$. Let $p_\phi: \tM_\phi
\rightarrow M$ be the covering introduced in section \ref{sec:topo}a. The
evaluation maps $\eta$ have natural lifts $\xi$,
\begin{equation} \label{eq:lift}
\xymatrix{
&&
{\tM_\phi} \ar[d]^{p_\phi} \\
{\MM_k(y_-,y_+,b,J)} \ar[urr]^-{\xi} \ar[rr]^-{\eta}
&&
M.
}
\end{equation}
Namely, $\xi(u)$ is the pair $(x,[c])$ consisting of $x = u(0,1)$ and the homotopy
class $[c]$ of the path $c(t) = u(0,t)$ from $\phi(x) = u(0,0)$ to $x$. In the same
way as before, one sees that $\xi$ extends naturally to a continuous map $\bar{\xi}$
on the Gromov-Floer compactification, which is a lift of $\bar{\eta}$.

Consider the pullback Morse function $\tilde{f} = f \circ p_\phi$ and metric on
$\tM_\phi$. Then $p_\phi^{-1}W^u(f,x_-)$ is the disjoint union of $W^u(\tilde{f},\tx_-)$
over all preimages $\tx_-$ of $x_-$. Together with the commutativity of \eqref{eq:lift}
this yields that
\[
\eta^{-1} W^u(f,x_-) = \bigsqcup_{p_\phi(\tx_-) = x_-} \xi^{-1}
W^u(\tilde{f},\tx_-) \subset \MM_k(y_-,y_+,b,J).
\]
In particular, if $k$ equals the Morse index of $x_-$ (and hence of $\tx_-$),
the sets $\xi^{-1}W^u(\tilde{f},\tx_-) \subset \MM_k(y_-,y_+,b,J)$ are finite, and
only a finite number of them is nonempty as one ranges over all $\tx_-$. Denoting
by $\tilde{q}(\tx_-,y_-,y_+)$ the number of points mod two in these sets, one defines
\begin{equation} \label{eq:lifted-chain}
\begin{split}
& CM^*(\tilde{f}) \otimes CF_*(\phi,b) \longrightarrow CF_*(\phi,b), \\ & \gen{\tx_-}
\otimes \gen{y_-} \longmapsto \sum_{y_+} \tilde{q}(\tx_-,y_-,y_+) \gen{y_+}.
\end{split}
\end{equation}
From looking at the one-dimensional spaces $\eta^{-1}W^u(\tilde{f},\tx_-)$ and applying
to their ends the same argument as before, one gets
\begin{multline*}
 \sum_{\tx} \tilde{l}(\tx_-,\tx) \tilde{q}(\tx,y_-,y_+) +
 \sum_{y} n(y_-,y) \tilde{q}(\tx_-,y,y_+) \\ +
\sum_{y} n(y,y_+) \tilde{q}(\tx_-,y_-,y) = 0 \qquad \qquad
\end{multline*}
with $\tilde{l}$ as defined in section \ref{sec:topo}c, so that $\tilde{\ast}$ is a
chain homomorphism. Moreover, by definition $\sum_{p_\phi(\tx_-) = x_-} \tilde{q}(\tx_-,
y_-,y_+) = q(x_-,y_-,y_+)$, so that the maps \eqref{eq:chain} and \eqref{eq:lifted-chain}
fit into a commutative diagram
\[
\xymatrix{
 {CM^*(\tilde{f}) \otimes CF_*(\phi,b)} \ar[r] &
 {CF_*(\phi,b)} \\
 {CM^*(f) \otimes CF_*(\phi,b)} \ar[u]^{p_\phi^* \otimes \id} \ar[ru]
}
\]
Therefore \eqref{eq:lifted-chain} induces a map $\tilde{\ast}: H^*(\tM_\phi)
\otimes HF_*(\phi) \rightarrow HF_*(\phi)$ which satisfies the conditions of
Lemma \ref{th:viterbo}. Having that, we can now prove the two theorems stated
at the beginning of the paper. If $g = [\phi]$ is nontrivial,
$H^2(\tM_\phi;\Z/2)$ is zero by Lemma \ref{th:pullback2}, which means that $a
\ast y = p_\phi^*(a) \tilde{\ast} y = 0$ for $a \in H^2(M;\Z/2)$ and $y \in
HF_*(g)$, as claimed in Theorem \ref{th:one}. In the same way, Lemma
\ref{th:viterbo} and Lemma \ref{th:pullback1} imply Theorem \ref{th:two}.


\providecommand{\bysame}{\leavevmode\hbox to3em{\hrulefill}\thinspace}

\end{document}